\documentstyle[12pt]{amsart}
\newtheorem{theorem}{Theorem}

\newtheorem{proposition}[theorem]{Proposition}

\newtheorem{corollary}[theorem]{Corollary}
\newtheorem{remar}[theorem]{Remark}
\newenvironment{proof}{Proof:\ \ \ }{\QED}
\newenvironment{remark}{\begin{remar}\rm}{\end{remar}}
\newcommand{\QED}{{\unskip\nobreak\hfil\penalty50%
\hskip1em\hbox{}\nobreak\hfil $\Box$%
\parfillskip=0pt \finalhyphendemerits=0 \par\medskip\noindent}}
\newcommand{\bfind}[1]{\index{#1}{\bf #1}}
\newcommand{\gloss}[1]{#1\glossary{\protect #1}}

\newcommand{\pars}{\par\smallskip}
\newcommand{\ovl}[1]{\overline{#1}}
\newcommand{\pH}{\mathop{\rule[-1pt]{0pt}{1pt}\mbox{\large\bf H}}}
\newcommand{\subsetuneq}{\mathrel{\raisebox{.8ex}{\footnotesize%
$\displaystyle\mathop{\subset}_{\not=}$}}}
\newcommand{\isom}{\simeq}
\newcommand{\N}{{\Bbb N}}

\newcommand{\R}{{\Bbb R}}
\begin{document}
\title{Exponentiation in power series fields}
\author{Franz-Viktor Kuhlmann, Salma Kuhlmann and Saharon Shelah}
\date{17.\ 5.\ 1996}
\thanks{The second author was supported by a Deutsche
Forschungsgemeinschaft fellowship. The third author was partially
supported by the Edmund Landau Center for research in Mathematical
Analysis, and supported by the Minerva Foundation (Germany). Publication
number 601}
\address{Mathematisches Institut der Universit\"at Heidelberg\\
         Im Neuenheimer Feld 288\\
         D-69120 Heidelberg, Germany}
\email{fvk@@harmless.mathi.uni-heidelberg.de}
\address{Department of Mathematics\\
         The Hebrew University of Jerusalem\\
         Jerusalem, Israel}
\email{shelah@@sunrise.huji.ac.il}

\begin{abstract}
We prove that for no nontrivial ordered abelian group $G$, the ordered
power series field $\R((G))$ admits an exponential, i.e.\ an
isomorphism between its ordered additive group and its ordered
multiplicative group of positive elements, but that there is a
non-surjective logarithm. For an arbitrary ordered field $k$, no
exponential on $k((G))$ is compatible, that is, induces an exponential
on $k$ through the residue map. This is proved by showing that certain
functional equations for lexicographic powers of ordered sets are not
solvable.
\end{abstract}
\maketitle

%
\section{Introduction}
Let $k$ be an ordered field and $G$ a nontrivial ordered abelian group,
then the (``generalized'') power series field $K=k((G))$ admits at least
one nonarchimedean order. Further, $K$ is real closed if and only if $k$
is real closed and $G$ is divisible. This provides a very simple and
elegant method of constructing nonarchimedean ordered real closed
fields. On the other hand, power series fields were already studied
by Levi-Civita [LC], later also by H.~Hahn [H], A.~Robinson [R] and
many others, in an attempt to develop function theory over
nonarchimedean fields. One of the first concerns was to define
elementary functions (e.g.\ the logarithm $\log$, or equivalently, its
inverse $\exp$) on those fields. It was already known to Levi-Civita
that $\log$ is definable through its Taylor expansion on the positive
units of $\R ((G))$, for archimedean ordered $G$ (cf.\ the discussion in
[L1], [L2]). This was generalized to arbitrary $G$ by B.~H.~Neumann [N].
But the problem of defining a logarithm from the group $K^{>0}$ of
positive elements onto $K$ remained open. We answer this problem to the
negative (Theorem \ref{nonexg}). In fact, we show that the domain of the
logarithm can be extended to $K^{>0}$ if $G$ is divisible, but that a
logarithm on $K^{>0}$ will never be surjective onto $K$.

For an \bfind{exponential} $f$ on the ordered field $(K,<)$ we only
require $f$ to be an isomorphism between its ordered additive group
$(K,+,0,<)$ and its ordered multiplicative group $(K^{>0},\cdot,1,<)$ of
positive elements. If $K=k((G))$, then we say that $f$ is
\bfind{compatible} if it induces an exponential on $k$ through the
canonical residue map (see Section 2 for details). We shall prove:
\begin{theorem}                             \label{nonexg}
Let $k$ be an ordered field and $G$ a nontrivial ordered abelian group.
Let $<$ be any order on $K=k((G))$. Then $(K,<)$ does not admit any
compatible exponential. If $(k,<)$ is archimedean, then $(K,<)$ admits
no exponential at all.
\end{theorem}
Theorem \ref{nonexg} shows that the construction method for real closed
fields described above is not available for exponential fields. Note
that there is an exponential on the surreal numbers (cf.\ [G]), but
this ``power series field'' is a proper class. For an even stronger
version of Theorem~\ref{nonexg}, see Theorem~\ref{strongerthm}.

The key to our result is the fact that every group complement of the
valuation ring in $K=k((G))$ is a lexicographic product of ordered
abelian groups. Let us recall the definition of lexicographic products.
Let $\Gamma$ and $\Delta_\gamma\,$, $\gamma\in\Gamma$ be totally ordered
sets. For every $\gamma\in\Gamma$, we fix a distinguished element $0\in
\Delta_\gamma\,$. The \bfind{support} of $a=(\delta_\gamma)_{\gamma\in
\Gamma}\in\prod_{\gamma\in\Gamma}\Delta_\gamma$, denoted by
\gloss{$\mbox{\rm supp} (a)$}, is the set of all $\gamma\in\Gamma$ for
which $\delta_\gamma\ne 0$. As a set, we define $\pH_{\gamma\in\Gamma}
\Delta_\gamma$ to consist of all $(\delta_\gamma)_{\gamma\in\Gamma}$
with well ordered support. The \bfind{lexicographic order} on
$\pH_{\gamma\in\Gamma}\Delta_\gamma$ is introduced as follows. Given
$a$ and $b=(\delta'_\gamma)_{\gamma\in\Gamma}\in\pH_{\gamma\in\Gamma}
\Delta_\gamma$, observe that $\mbox{\rm supp}(a)\cup\mbox{\rm supp}(b)$
is well ordered. Let $\gamma_0$ be the least of all elements
$\gamma\in\mbox{\rm supp}(a) \cup \mbox{\rm supp}(b)$ for which
$\delta_\gamma\ne\delta'_\gamma\,$. We set $a<b\,:\Leftrightarrow\,
\delta_{\gamma_0}<\delta'_{\gamma_0}\,$. Then $(\pH_{\gamma\in\Gamma}
\Delta_\gamma,<)$ is a totally ordered set, the \bfind{lexicographic
product} of the ordered sets $\Delta_\gamma\,$. If $\Delta_\gamma=
\Delta$ for all $\gamma\in\Gamma$ then we write $\Delta^\Gamma$ for
their lexicographic product; it consists of all maps from $\Gamma$ to
$\Delta$ with well ordered support.

If all $\Delta_\gamma$ are totally ordered abelian groups, then we can
take the distinguished elements $0$ to be the neutral elements of the
groups $\Delta_\gamma\,$. Defining addition on $\pH_{\gamma\in\Gamma}
\Delta_\gamma$ componentwise, we obtain a totally ordered abelian group
$(\pH_{\gamma\in\Gamma}\Delta_\gamma,+,0,<)$, the \bfind{Hahn product}
of the ordered groups $\Delta_\gamma\,$.

\pars
In Section~\ref{sectlexp}, we prove the following theorem and
explain how it relates to the surjectivity of a logarithm.
\begin{theorem}                                     \label{shelah}
Let $\Gamma$ and $\Delta$ be totally ordered sets without greatest
element, and fix an element $0\in\Delta$. Suppose that
$\Gamma'$ is a cofinal subset of $\Gamma$ and that $\iota\colon
\Gamma'\,\rightarrow\,\Delta^\Gamma$ is an order preserving embedding.
Then the image $\iota\Gamma'$ is not convex in $\Delta^\Gamma$.
\end{theorem}
The same holds for an order preserving embedding $\iota\colon\Gamma'
\,\rightarrow\, \pH_{\gamma\in\Gamma}\Delta_\gamma$ and already under
the condition that $\Gamma$ has no greatest element and $0$ is not the
greatest element of $\Delta_\gamma$ for any $\gamma\in\Gamma$ (cf.\
[K--K--S]). If we drop the condition that $\Gamma$ has no greatest
element, the situation changes drastically. Suitably chosen ordered sets
$\Gamma$ and $\Delta$ will even admit an isomorphism $\Gamma\isom
\Delta^\Gamma$. We study this situation and related questions in
[K--K--S].

%
%
\section{Preliminaries on left logarithms}  \label{sectprel}
Let $G$ be a totally ordered abelian group. The set of archimedean
classes $[a]$ of nonzero elements $a\in G$ is totally ordered by setting
$[a]<[b]$ if $|a|\gg |b|$. The map $v:\>a\mapsto [a]$ is called the
\bfind{natural valuation} of $G$. It satisfies the triangle inequality
$v(a-b)\geq\min\{va,vb\}$ and $v(-a)=va$ as well as
\begin{equation}                            \label{<=v<=}
a\leq b\leq 0\>\vee\> a\geq b\geq 0\;\;\Rightarrow\;\; va\leq vb\;.
\end{equation}
In this paper, $(K,<)$ will always be a totally ordered
field. We let $v$ denote the natural valuation on its additive group
$(K,+,0,<)$. In this case, $vK:=v(K\setminus\{0\})$ forms a totally
ordered abelian group endowed with the addition $[a]+[b]:=[ab]$, and
$v$ is a field valuation. For more information on natural valuations,
see [K].

Let $w$ be any field valuation on $K$. The \bfind{value group} of
$(K,w)$ will be denoted by $wK$ and its \bfind{residue field} by $Kw$.
Further, $w$ is \bfind{convex with respect to $<$} if it satisfies
(\ref{<=v<=}). The \bfind{valuation ring} $R_w=\{a\in K\mid wa\geq 0\}$
of a convex valuation $w$ is convex in $K$, and so is its
\bfind{valuation ideal} $I_w=\{a\in K\mid wa>0\}$. Further, the set
${\cal U}_w^{>0}:=\{a\in K\mid wa=0\,\wedge\, a>0\}$ of \bfind{positive
units} of $R_w$ is a convex subgroup of $(K^{>0},\cdot,1,<)\,$. Note
that $w$ is convex if and only if $R_v\subseteq R_w$ (i.e., $w$ is a
\bfind{coarsening} of $v$), in which case its value group $wK$ is the
quotient of $vK$ by a convex subgroup.

If $K$ admits an exponential, then its multiplicative group of positive
elements is divisible (since the additive is). In order to prove
Theorem~\ref{nonexg}, we can thus always assume divisibility. As in
[K] (Lemma~3.4 and Theorem~3.8), we then have the following
representations as lexicographic sums:
\begin{equation}                            \label{A}
(K,+,0,<)\>\isom\> {\bf A}_w\,\amalg\, (R_w,+,0,<)
\end{equation}
where ${\bf A}_w$ is an arbitrary group complement of $R_w$ in $(K,+)$,
and analogously,
\begin{equation}                            \label{B}
(K^{>0},\cdot ,1,<)\>\isom\> {\bf B}_w\,\amalg\,
({\cal U}_w^{>0},\cdot ,1,<)
\end{equation}
where ${\bf B}_w$ is an arbitrary group complement of ${\cal U}_w^{>0}$
in $(K^{>0},\cdot\,)\,$. Endowed with the restriction of the ordering,
${\bf A}_w$ and ${\bf B}_w$ are unique up to isomorphism. In view of
(\ref{<=v<=}) and the fact that $w(-a)=wa$, the map
\begin{equation}                            \label{-v}
-w:\>(K^{>0},\cdot,1,<)\;\rightarrow\; (wK,+,0,<)\,,\hspace{1cm}
-wa = wa^{-1}
\end{equation}
is a surjective group homomorphism preserving $\leq\,$, with kernel
${\cal U}_w^{>0}$. We find that every complement ${\bf B}_w$ is
isomorphic to $(wK,+,0,<)$ through the map $-w$.

An exponential $f$ on $K$ will be called \bfind{compatible with $w$} if
it satisfies that
$f(R_w)= {\cal U}_w^{>0}$ and $f(I_w)= 1+I_w\,$. Since $Kw=R_w/I_w$ and
$(Kw^{>0},\cdot ,1,<)$ $=({\cal U}_w^{>0},\cdot ,1,<)/\,1+I_w\,$, this
means that $f$ induces canonically an exponential $fw:\>(Kw,+,0,<)
\rightarrow (Kw^{>0},\cdot ,1,<)$. The canonical valuation $w$ of a
power series field $K=k((G))$ has value group $G$ and residue field $k$.
Further, it is henselian. Consequently, $w$ is convex with respect to
every order $<$ on $K$ (cf.\ [KN--WR]). Hence, $f$ is compatible on
$k((G))$ if and only if it is compatible with $w$.

\begin{remark}
If an ordered field $K$ admits an exponential, then it admits an
exponential which is compatible with the natural valuation (cf.\ [K],
Section 3.3).
\end{remark}

Let $w$ be any convex valuation on $K$. Every compatible exponential $f$
decomposes into two isomorphisms of ordered groups:
\[\begin{array}{rrcl}
f_R:& (R_w,+,0,<) & \rightarrow & ({\cal U}_w^{>0},\cdot ,1,<)\\
f_L:& {\bf A}_w & \rightarrow & {\bf B}_w\;.
\end{array}\]
Conversely, in view of (\ref{A}) and (\ref{B}), such isomorphisms $f_R$
and $f_L$ can be put together to obtain an exponential compatible with
$w$. The inverse $f_L^{-1}$ is called a \bfind{left logarithm}, and
$f_R^{-1}$ a \bfind{right logarithm}. Through the isomorphism
(\ref{-v}), every isomorphism
\[h:\;(wK,+,0,<)\,\rightarrow\,{\bf A}_w\]
gives rise to a left logarithm $h\circ -w$. Conversely, given a left
logarithm $f_L^{-1}$, the map $f_L^{-1} \circ (-w)^{-1}$ is such an
isomorphism $h$. This correspondence
motivates the following definition: a \bfind{logarithmic cross-section}
of an ordered field $(K,<)$ with respect to a convex valuation $w$ is an
order preserving embedding $h$ of $wK$ into an additive group complement
of the valuation ring (that is, an embedding $h$ of $wK$ into the
additive group $(K,+,0,<)$ satisfying $wh(g) < 0$ for all $g\in wK$).
Thus, every left logarithm induces a logarithmic cross-section which is
surjective (i.e., $h(wK)$ is an additive group complement to the
valuation ring), and vice-versa.

Our goal in the next section is to show that power series fields
{\it always} admit logarithmic cross-sections, but {\it never}
surjective ones.

%
%
\section{Lexicographic products and logarithmic
cross-sections}                                     \label{sectlexp}
{\bf Proof of Theorem~\ref{shelah}:}\ \
Assume that $\Gamma$ and $\Delta$ are totally ordered sets, and fix an
element $0\in\Delta$. Assume further that $\Delta$ has no greatest
element, so that we can choose a map $\tau\colon\Delta\,\rightarrow\,
\Delta$ such that $\tau\delta>\delta$ for all $\delta\in\Delta$.
For every well ordered set $S\subset \Gamma$ and every
$d=(d_\gamma)_{\gamma\in\Gamma} \in\Delta^\Gamma$, set
\[d\oplus S=(d'_\gamma)_{\gamma\in\Gamma}\;\mbox{ where }
d'_\gamma:=\left\{\begin{array}{cl} d_\gamma &\mbox{if }\gamma\notin S\\
\tau d_\gamma & \mbox{if } \gamma\in S\;. \end{array}\right.\]
Observe that the support of $d\oplus S$ is contained in
$\mbox{\rm supp}(d)\cup S$ and thus, it is again well ordered.
Further, if $S,S'\subset\Gamma$ are well ordered sets (or empty), then
\begin{equation}                            \label{oplusS}
S\subsetuneq S'\;\Rightarrow\;d\oplus S<d\oplus S'\;.
\end{equation}

\pars
Now suppose that $\Gamma$ has no greatest element, $\Gamma'$ is a
cofinal subset of $\Gamma$ and $\iota\colon\Gamma'\,\rightarrow\,
\Delta^\Gamma$ is an order preserving embedding such that the image
$\iota\Gamma'$ is convex in $\Delta^\Gamma$. We wish to deduce a
contradiction.

By induction on $n\in\N$, we define elements $\gamma_0^{(n)} \in
\Gamma'$. We choose an arbitrary $\gamma_0^{(1)}\in \Gamma'$. Having
already constructed $\gamma_0^{(n)}$, we carry through the following
induction step. Since $\Gamma$ has no greatest element, the same holds
for $\Gamma'$, and there is some $\alpha^{(n)}\in\Gamma'$ such that
$\gamma_0^{(n)}<\alpha^{(n)}$. Hence, $\iota\gamma_0^{(n)}<\iota
\alpha^{(n)}$. Let $\beta^{(n)}\in \Gamma$ be the least element of
$\mbox{\rm supp}(\iota\gamma_0^{(n)}) \cup\mbox{\rm supp}
(\iota\alpha^{(n)})$ for which
\[(\iota\gamma_0^{(n)})_{\beta^{(n)}}\> < \>
(\iota\alpha^{(n)})_{\beta^{(n)}}\;.\]
Since $\Gamma$ has no greatest element and $\Gamma'$ is a cofinal
subset, we can choose $\gamma_0^{(n+1)} \in \Gamma'$ such that
$\beta^{(n)}<\gamma_0^{(n+1)}$.

If $S\subset \Gamma$ is a well
ordered set with least element $\gamma_0^{(n+1)}$, then
\begin{equation}                            \label{<oplus}
\iota\gamma_0^{(n)}\><\>\iota\gamma_0^{(n)}\oplus S\><\>
\iota\alpha^{(n)}\;.
\end{equation}
Indeed, $(\iota\gamma_0^{(n)}\oplus S)_\beta=
(\iota\gamma_0^{(n)})_\beta$ for every $\beta<\gamma_0^{(n+1)}$. In
particular,
\[(\iota\gamma_0^{(n)}\oplus S)_{\beta^{(n)}}=
(\iota\gamma_0^{(n)})_{\beta^{(n)}}
<(\iota\alpha^{(n)})_{\beta^{(n)}}\;,\]
which implies the second inequality of (\ref{<oplus}).
Its first inequality follows from (\ref{oplusS}).

The image of $\Gamma'$ in $\Delta^\Gamma$ being convex, (\ref{<oplus})
yields that also $\iota\gamma_0^{(n)}\oplus S$ lies in this image. Thus,
$\iota^{-1}(\iota\gamma_0^{(n)}\oplus S)$ is a well defined element of
$\Gamma'$.

Suppose now that for some ordinal number $\mu\geq 1$ we have chosen
elements $\gamma_\nu^{(n)}\in\Gamma'$, $\nu<\mu$, $n\in\N$, such that
for every fixed $n$, the sequence $(\gamma_\nu^{(n)})_{\nu<\mu}$ is
strictly increasing. Then we set
\[\gamma_\mu^{(n)}:=\iota^{-1}(\iota\gamma_0^{(n)}\oplus
\{\gamma_\nu^{(n+1)} \mid\nu<\mu\})\in\Gamma'\]
for every $n\in\N$. If $\lambda<\mu$, then $\{\gamma_\nu^{(n+1)} \mid
\nu<\lambda\}\subsetuneq\{\gamma_\nu^{(n+1)} \mid\nu<\mu\}$ and thus,
$\gamma_\lambda^{(n)}<\gamma_\mu^{(n)}$ by (\ref{oplusS}). So for every
ordinal number $\mu$, the sequences $(\gamma_\nu^{(n)})_{\nu<\mu}$ can
be extended. We obtain strictly increasing sequences of arbitrary
length, contradicting the fact that their length is bounded by the
cardinality of $\Gamma'$.                                         \QED

Now we apply Theorem~\ref{shelah} to logarithmic cross-sections of the
power series field $K=k((G))$ with canonical valuation $w$. One of the
complements for the valuation ring $R_w=k[[G]]$ is the Hahn product
$\pH_{G^{<0}} (k,+,0,<)$, which we will denote by $k^{G^{<0}}$. Since the
complements are unique up to isomorphism, a surjective logarithmic
cross-section $h$ with respect to $w$ would induce an isomorphism
$G\simeq k^{G^{<0}}$. This in turn would imply that $G^{<0}$ has no
greatest element and would give rise to an embedding of
$G^{<0}$ in $k^{G^{<0}}$ with convex image, which contradicts
Theorem~\ref{shelah}. So we have proved:
\begin{theorem}                             \label{nosurj}
Let $k$ be an ordered field and $G$ a nontrivial ordered abelian group.
Further, let $w$ be the canonical valuation on $K=k((G))$ and $<$ any
order on $K$. Then $(K,<)$ admits no surjective logarithmic
cross-section with respect to $w$.
\end{theorem}
This theorem implies Theorem \ref{nonexg}. Indeed, a compatible
exponential of $K$ would induce a surjective logarithmic cross-section
with respect to $w$, which is impossible. If $(k,<)$ is archimedean,
then $w$ will coincide with the natural valuation $v$ of $(k((G)),<)$.
So the second assertion of Theorem \ref{nonexg} follows by Remark~3.

\pars
If $G$ is an ordered abelian group, then we denote its natural valuation
by $v_G\,$. For the definition of the archimedean components $B_\gamma$
of $G$ (where $\gamma\in v_G G$), see [FU]. They are archimedean ordered
abelian groups. Hahn's embedding theorem states that there is an order
preserving group embedding $\rho$ of $G$ in the Hahn product
$\pH_{\gamma\in v_G G} B_\gamma\,$, if $G$ is divisible (cf.\ [H] or
[FU], IV, Theorem~16).
\begin{proposition}
Let $G$ be a nontrivial divisible ordered abelian group. Then $\R((G))$
admits a logarithmic cross-section. If every archimedean component of
$G$ embeds in the ordered additive group of $k$, then $k((G))$ admits a
logarithmic cross-section.
\end{proposition}
\begin{proof}
By taking representatives, we obtain an embedding $\sigma:\; v_G G
\rightarrow\; G^{<0}$; it is order preserving by (\ref{<=v<=}). Now
$\sigma$ lifts to an embedding $\hat\sigma:\; k^{v_G G} \rightarrow\;
k^{G^{<0}}$. If every archimedean component $B_\gamma$ of $G$ embeds in
$(k,+,0,<)$, then there is an embedding $\tau:\;\pH_{\gamma\in v_G G}
B_\gamma \rightarrow\;\pH_{v_G G} (k,+,0,<)=k^{v_G G}$. So
$h=\hat\sigma\circ\tau \circ\rho$ is the required logarithmic
cross-section for $k((G))$. Since every archimedean ordered abelian
group embeds in $\R$, the first assertion follows from the second.
\end{proof}
\begin{corollary}
If $G$ is nontrivial and divisible, then the real closed field
$K=\R((G))$ admits a non-surjective logarithm, i.e., an embedding
$(K^{>0},\cdot,1,<)\rightarrow (K,+,0,<)$.
\end{corollary}
For the proof, note that a right logarithm on $\R((G))$ always exists:
it is defined on the positive units of the valuation ring $\R[[G]]$
through the logarithmic power series (cf.\ [A]). In combination with a
non-surjective logarithmic cross-section, it gives rise to
the desired (non-surjective) logarithm. By taking the union over a
suitable countable ascending chain of power series fields with
non-surjective logarithms, we can obtain a surjective logarithm. Using
this construction, we prove in [K--K2] the existence of exponential
fields with arbitrary given exponential rank ($=$ the order type of the
set of all convex valuations compatible with the exponential).

\begin{remark}
\small\rm
Kaplansky [KA] has shown that a valued field is maximal (i.e., admits
no proper immediate extensions) if and only if every pseudo Cauchy
sequence admits a limit. The same principle was proved by Fleischer
[F] for valued abelian groups. It can also be proved for certain
classes of valued modules. At first sight, one might believe that this
principle holds for all (reasonable) valued structures. But the
nonarchimedean exponential fields with their natural valuation
constitute a counterexample:

There are maximal naturally valued exponential fields (i.e., they do not
admit proper immediate extensions to which also the exponential
extends). These are precisely the
exponential fields whose natural valuation $v$ is complete: On the one
hand, it was remarked in [K] that if $(L,v)\supset (K,v)$ is immediate
and the exponential extends from $K$ to $L$, then $(K,v)$ is dense in
$(L,v)$. On the other hand, if $(K,v)$ is dense in $(L,v)$, then an
exponential of $K$ extends to $L$ by continuity. Hence, the completion
of a nonarchimedean exponential field with respect to its natural
valuation is the maximal immediate extension as a naturally valued
exponential field. But by our nonexistence result, it cannot be a power
series field. On the other hand, Kaplansky has also shown in [KA] that
a valued field $(K,w)$ of residue characteristic 0 is a power series
field with canonical valuation $w$ if and only if it is maximal. (Note
that the natural valuation has residue characteristic 0 since the
residue field is ordered.) Hence, a maximal naturally valued
exponential field is not maximal as a valued field.

An argument similar to that used in establishing Theorem~\ref{nosurj}
shows that a Hahn group (i.e.\ a maximally valued group) cannot be an
exponential group in the sense of [K], and thus cannot be the
natural value group of an exponential field. Further consequences of
Theorem~\ref{shelah} for exponential groups will be studied in
a subsequent paper.
\end{remark}

Under the hypothesis of Theorem~\ref{nonexg} we can prove that an
exponential cannot even be compatible with any nontrivial coarsening
$w'$ of $w$: Since $K$ is a power series field with canonical valuation
$w$, it can also be written as a power series field $(Kw)((w'K))$ with
canonical valuation $w'$, and from Theorem~\ref{nonexg} it follows that
no exponential can be compatible with $w'$. We have seen in the above
remark that we can talk about maximal valuations instead of power series
fields. So we can restate our result as follows:\ \ {\it If the ordered
field $K$ admits an exponential $f$, then there is no nontrivial
coarsening of its natural valuation $v$ which is maximal and compatible
with $f$.} We prove the following generalization:
\begin{theorem}                             \label{strongerthm}
Let $f$ be an exponential on the ordered field $K$ and $w$ a coarsening
of the natural valuation $v$ of $K$ such that $f$ is compatible with
$w$. Then there is no coarsening $\tilde{w}$ of $w$ such that the
valuation $\ovl{w}=w/\tilde{w}$ induced by $w$ on the residue field
$K\tilde{w}$ is nontrivial and $(K\tilde{w},\ovl{w})$ is maximal.
\end{theorem}
\begin{proof}
Suppose to the contrary that such a coarsening $\tilde{w}$ exists. We
have that $R_w \subset R_{\tilde{w}}$. Let $\ovl{\bf A}$ be a group
complement of $R_w$ in $R_{\tilde{w}}$ and $\tilde{\bf A}$ a group
complement of $R_{\tilde{w}}$ in $(K,+,0,<)$. Then $\tilde{\bf A}\amalg
\ovl{\bf A}$ is a group complement of $R_w$ in $(K,+,0,<)$. Further, $f$
induces an isomorphism $h$ from $G =wK$ onto $\tilde{\bf A}\amalg
\ovl{\bf A}$ as ordered groups. In particular, $G^{<0}$ has no greatest
element.

The value group of $\ovl{w}$ is isomorphic to a nontrivial convex
subgroup $\ovl{G}$ of $G$. Since $(K\tilde{w},\ovl{w})$ is maximal and
has residue field $(K\tilde{w}) w/\tilde{w}=Kw$, it is isomorphic to the
power series field $(Kw)((\ovl{G}))$. Hence, $\ovl{\bf A}$ is isomorphic
to a Hahn product $(Kw)^{\ovl{G}^{<0}}$. This yields an embedding of the
nontrivial convex subgroup $H:=\ovl{G}\cap h^{-1}(\ovl{\bf A})$ of
$\ovl{G}$ in $(Kw)^{\ovl{G}^{<0}}$. Under this embedding, the image of
the final segment $H^{<0}$ of $\ovl{G}^{<0}$ is convex in
$(Kw)^{\ovl{G}^{<0}}$. But $\ovl{G}^{<0}$ is a final segment of $G^{<0}$
and thus has no greatest element. This contradicts Theorem~\ref{shelah}.
\end{proof}

\end{document}